\documentclass[leqno,11pt]{amsart}
\usepackage{amsmath,amscd,amsthm,amsxtra}
\usepackage{epsfig,graphics,color,colortbl}
\usepackage{amssymb,latexsym}
\usepackage{mathabx}
\usepackage{mathrsfs,eucal,upgreek}
\usepackage[poly,all]{xy}
\usepackage{hyperref}
\usepackage{tikz-cd}
\usepackage{arydshln}

\newtheorem{thm}{\bf Theorem}[section]

\newtheorem{prop}[thm]{\bf Proposition}
\newtheorem{cor}[thm]{\bf Corollary}
\newtheorem{lem}[thm]{\bf Lemma}
\newtheorem{rem}[thm]{\bf Remark}

\numberwithin{equation}{section}

\newcommand{\bs}{\boldsymbol}

\newcommand{\B}{\mathbf{B}}

\newcommand{\W}{\mathcal{W}}
\newcommand{\cP}{\mathscr{P}}

\newcommand{\pf}{\noindent{\bfseries Proof. }}
\newcommand{\ov}{\overline}

\newcommand{\ba}{\bs{\rm a}}
\newcommand{\bb}{\bs{\rm b}}

\newcommand{\M}{{\mathcal{M}}}

\newcommand{\gl}{\mathfrak{gl}}
\newcommand{\Z}{\mathbb{Z}}
\newcommand{\C}{\mathbb{C}}
\newcommand{\h}{\mathfrak{h}}
\newcommand{\te}{\widetilde{e}}
\newcommand{\tf}{\widetilde{f}}
\newcommand{\g}{\mathfrak{g}}
\newcommand{\td}{\widetilde}

\newcommand{\mc}{\mathcal}
\newcommand{\mf}{\mathfrak}

\newcommand{\I}{\mathcal{I}}

\newcommand{\la}{\lambda}

\usepackage{tikz}
\usetikzlibrary{matrix,calc}

\begin{document}
\title[]
%{A simple proof of Lascoux's non-symmetric Cauchy identity}
{Lakshmibai-Seshadri paths and non-symmetric Cauchy identity}

\author{SEUNG-IL CHOI$^\dagger$}
\address{Department of Mathematical Sciences, Seoul National University, Seoul 08826, Korea}
\email{ignatioschoi@snu.ac.kr}
\thanks{$^\dagger$Partially supported by Samsung Science and Technology Foundation under Project Number SSTF-BA1501-01.}

\author{JAE-HOON KWON$^{\dagger\dagger}$}
\address{Department of Mathematical Sciences, Seoul National University, Seoul 08826, Korea}
\email{jaehoonkw@snu.ac.kr}
\thanks{$^{\dagger\dagger}$Partially supported by Research Resettlement Fund for the new faculty of Seoul National University.}


\begin{abstract}
We give a simple crystal theoretic interpretation of the Lascoux's expansion of a non-symmetric Cauchy kernel $\prod_{i+ j\leq n+1}(1-x_iy_j)^{-1}$, which is given in terms of Demazure characters and atoms. We give a bijective proof of the non-symmetric Cauchy identity using the crystal of Lakshmibai-Seshadri paths, and extend it to the case of continuous crystals.
\end{abstract}

\maketitle
%\setcounter{tocdepth}{1}
%\tableofcontents

\section{Introduction}
Let $S(\C^n\otimes \C^n)$ be the symmetric algebra over the complex numbers generated by $\C^n\otimes\C^n$.
It is well-known that $S(\C^n\otimes \C^n)$ admits a natural action of $GL_n \times GL_n$, and the following multiplicity-free decomposition:
\begin{equation}\label{eq:Cauchy decomp}
S(\C^n\otimes \C^n) \cong \bigoplus_{\la}V(\la)\otimes V(\la),
\end{equation}
where the sum is over the partitions $\la$ with length no more than $n$, and $V(\la)$ denotes the irreducible polynomial $GL_n$-module corresponding to $\la$. This decomposition implies the Cauchy identity, where the Cauchy kernel $\prod_{1\leq i,j\leq n}(1-x_iy_j)^{-1}$ is given as a sum of products of two Schur polynomials in $x=\{\,x_1,\ldots,x_n\,\}$ and $y=\{\,y_1,\ldots,y_n\,\}$ corresponding to $\la$. Similarly, the exterior algebra generated by $\C^n\otimes \C^n$ yields the dual Cauchy identity, which expands the inverse of the Cauchy kernel after replacing $y$ with $-y$.
%\begin{equation*}
%\frac{1}{\prod_{1\leq i,j\leq n}(1-x_iy_j)}=\sum_{\lambda}s_\lambda(x)s_\lambda(y)
%\end{equation*}
%where the sum is over the partitions with length no more than $n$, and $s_\la(x)$ and $s_\la(y)$ are the Schur polynomials in $x_i$'s and $y_j$'s respectively.

There is a non-symmetric analogue of the Cauchy identity introduced in \cite{FuLa,La}, which
expands ${\prod_{i+ j\leq n+1}(1-x_iy_j)}^{-1}$ as follows:
\begin{equation}\label{eq:non-sym Cauchy kernel}
\begin{split}
\frac{1}{\prod_{i+ j\leq n+1}(1-x_iy_j)}
&=\sum_{\la\in\cP_n}\sum_{w\in W/W_\la}
{K}_{w\la}(x)\widehat{K}_{ww_0\la}(y),
\end{split}
\end{equation}
where $K_{w\la}(x)$ and $\widehat{K}_{ww_0\la}(y)$ are the Demazure characters and Demzure atoms in $x$ and $y$, respectively (see Section \ref{sec:Cauchy} for more details).
Like the usual Cauchy identity, one may regard \eqref{eq:non-sym Cauchy kernel} as a dual of the expansion of the product $\prod_{i+j\leq n+1}(x_i+y_j)$ given as a sum of products of two Schubert polynomials in $x$ and $y$.

%The purpose of this article is to understand a representation theoretic meaning of the non-symmetric Cauchy identity \eqref{eq:non-sym Cauchy kernel}, the question on which was raised by Lascoux \cite{La} when he introduced the identity.

%The proof is based on the theory of crystal base (cf.~\cite{Kas95}) and the combinatorics of Littelmann path model or the the crystal of Lakshmibai-Seshadri paths \cite{Li95}.

In this paper, we give a simple crystal theoretic interpretation and bijective proof of the non-symmetric Cauchy identity \eqref{eq:non-sym Cauchy kernel}.

Let us recall that the crystal associated to the symmetric algebra $S(\C^n\otimes \C^n)$ or the quantized coordinate ring of $n\times n$ matrices can be realized as the set $\M$ of $n\times n$ non-negative integral matrices, which is a $(\gl_n,\gl_n)$-bicrystal \cite{DK,La,vL}. The main idea is to regard $\M$ as a $\gl_n$-crystal via the tensor product rule, which corresponds to the diagonal action of $GL_n$ on $S(\C^n\otimes\C^n)$ in \eqref{eq:Cauchy decomp}. We consider the $\gl_n$-subcrystal of $\M$ generated by the matrices corresponding to the $(GL_n,GL_n)$-highest weight vectors, which is given by the set $\M^{low}$ of lower triangular matrices in $\M$.
We show that there exists a bijection from $\M^{low}$ to the set of pairs of Lakshmibai-Seshadri paths lying in the two-fold tensor power of a highest weight crystal and satisfying certain ordering conditions (Theorem \ref{thm:main}). Then the bijection recovers the identity \eqref{eq:non-sym Cauchy kernel} if we take the ``refined characters" of both sides of the bijection by distinguishing the weights from each tensor factor in terms of $x$ and $y$ (Corollary \ref{eq:non sym Cauchy id-2}). In this sense, it is interesting to note that \eqref{eq:non-sym Cauchy kernel} is a refinement of the classical Littlewood identity \eqref{eq:Littlewood identity}.

We also see from our proof that the right-hand side of \eqref{eq:non-sym Cauchy kernel} appears naturally in a more general setting. Indeed, if we consider the quantum group $U_q(\g)$ and its quantum coordinate ring $A_q(\g)$ associated to a semisimple Lie algebra $\g$ \cite{Kas93-1}, then the refined character (in the above sense) of the $U_q(\g)$-module generated by the $(U_q(\g)\otimes U_q(\g))$-highest weight vectors via comultiplication is given by a sum of the products of opposite Demazure characters and Demazure atoms for $\g$ (Remark \ref{rem:semisimple case}).
Finally, as an application of our bijective proof we give an analogue of the non-symmetric Cauchy decomposition for continuous crystals  \cite{BBO,Kas02} (Theorem \ref{thm:main-2}).

We remark that the bijective proof of \eqref{eq:non-sym Cauchy kernel} in this paper is different from those in \cite{AE,La}, though the theory of crystals (of tableaux) is still used there.
We hope that it would be helpful for a more direct representation theoretic interpretation of the non-symmetric Cauchy identity, the question on which was raised by Lascoux \cite{La} when he introduced it.

\section{Demazure crystals and path model}\label{subsec:crystal}
\subsection{Crystals}%Review on Cartan data for KM alg, (Demazure/opposite) Crystals
Let us give a brief review on crystals (cf. \cite{HK,Kas95}). We denote by $\Z_+$ and $\mathbb{R}_+$ the sets of non-negative integers and non-negative real numbers, respectively.
Let $\g$ be the Kac-Moody algebra associated to a symmetrizable generalized Cartan matrix $A =(a_{ij})_{i,j\in I}$ for some finite index set $I$.
Let $P^\vee$ be the dual weight lattice,
$P = {\rm Hom}_\Z( P^\vee,\Z)$ the weight lattice,
$\Pi^\vee=\{\,h_i\,|\,i\in I\,\}$ the set of simple coroots, and
$\Pi=\{\,\alpha_i\,|\,i\in I\,\}$ the set of simple roots of $\g$ such that $\langle \alpha_j,h_i \rangle=a_{ij}$ for $i,j\in I$. Let $P_+=\{\,\la\in \h^*\,|\,\langle \la,h_i \rangle\in \Z_+\,\}$ the set of dominant integral weights. Let $\h=\C\otimes_{\Z} P^\vee$ be the Cartan subalgebra of $\g$, and let $W\subset GL(\h^\ast)$ the Weyl group of $\g$ generated by $s_i$ for $i\in I$, where $s_i(\la)=\la-\la(h_i)\alpha_i$ for $\la\in \h^\ast$.

A crystal is a set
$B$ together with the maps ${\rm wt} : B \rightarrow P$,
$\varepsilon_i, \varphi_i: B \rightarrow \mathbb{Z}\cup\{-\infty\}$ and
$\te_i, \tf_i: B \rightarrow B\cup\{{\bf 0}\}$ for $i\in I$ satisfying the following:  for $b\in B$ and $i\in I$,
\begin{itemize}
\item[(1)]
$\varphi_i(b) =\langle {\rm wt}(b),h_i \rangle +
\varepsilon_i(b)$,

\item[(2)] $\varepsilon_i(\te_i b) = \varepsilon_i(b) - 1,\ \varphi_i(\te_i b) =
\varphi_i(b) + 1,\ {\rm wt}(\te_ib)={\rm wt}(b)+\alpha_i$ if $\te_i b \in B$,

\item[(3)] $\varepsilon_i(\tf_i b) = \varepsilon_i(b) + 1,\ \varphi_i(\tf_i b) =
\varphi_i(b) - 1,\ {\rm wt}({\tf_i}b)={\rm wt}(b)-\alpha_i$ if $\tf_i b \in B$,

\item[(4)] $\tf_i b = b'$ if and only if $b = \te_i b'$ for $b' \in B$,

\item[(5)] $\te_ib=\tf_ib={\bf 0}$ when $\varphi_i(b)=-\infty$,
\end{itemize}
where ${\bf 0}$ is a formal symbol and $-\infty$ is the smallest
element in $\Z\cup\{-\infty\}$ such that $-\infty+n=-\infty$
for all $n\in\Z$.
For $\mu\in\mf h^*$, let
$B_\mu=\{\, b\in B\,  | \, {\rm wt}(b)=\mu\,\}$.
We have $B=\bigsqcup_{\mu}B_\mu$. When $B_\mu$ is finite for all $\mu$, we define the character of $B$ by ${\rm ch}B=\sum_{\mu\in \mf h^*}|B_\mu|e^\mu$,
where $e^\mu$ is a basis element of the group algebra $\mathbb{Q}[\mf h^\ast]$.
%Unless otherwise specified, a crystal means a $\g$-crystal throughout the paper for simplicity.
%A crystal becomes an $I$-colored oriented graph, where $b\stackrel{i}{\rightarrow}b'$ if and only if $b'=\tf_{i}b$ for $b, b'\in B$ and $i\in I$.

Let $B_1$ and $B_2$ be crystals.
A tensor product $B_1\otimes B_2$
is a crystal, which is defined to be $B_1\times B_2$  as a set with elements  denoted by
$b_1\otimes b_2$, where
\begin{equation*}%\label{eq:tensor product of crystals}
\begin{split}
{\rm wt}(b_1\otimes b_2)&={\rm wt}(b_1)+{\rm wt}(b_2), \\
\varepsilon_i(b_1\otimes b_2)&= {\rm
max}\{\varepsilon_i(b_1),\varepsilon_i(b_2)-\langle {\rm
wt}(b_1),h_i\rangle\}, \\
\varphi_i(b_1\otimes b_2)&= {\rm max}\{\varphi_i(b_1)+\langle {\rm
wt}(b_2),h_i\rangle,\varphi_i(b_2)\},
\\
{\te}_i(b_1\otimes b_2)&=
\begin{cases}
{\te}_i b_1 \otimes b_2, & \text{if $\varphi_i(b_1)\geq \varepsilon_i(b_2)$}, \\
b_1\otimes {\te}_i b_2, & \text{if
$\varphi_i(b_1)<\varepsilon_i(b_2)$},
\end{cases}\\
{\tf}_i(b_1\otimes b_2)&=
\begin{cases}
{\tf}_i b_1 \otimes b_2, & \text{if  $\varphi_i(b_1)>\varepsilon_i(b_2)$}, \\
b_1\otimes {\tf}_i b_2, & \text{if $\varphi_i(b_1)\leq
\varepsilon_i(b_2)$},
\end{cases}
\end{split}
\end{equation*}
\noindent for $i\in I$. Here we assume that ${\bf 0}\otimes
b_2=b_1\otimes {\bf 0}={\bf 0}$. Given $b_i \in B_i $ ($i=1,2$), we write $b_1 \equiv b_2$ if there is an isomorphism of crystals $C(b_1) \rightarrow C(b_2)$ mapping $b_1$ to $b_2$, where $C(b_i)$ denotes the connected component of $b_i$ in $B_i$.

For an indeterminate $q$, let $U_q(\g)$ be the quantized enveloping algebra of $\g$ over $\mathbb{Q}(q)$. For $\la\in P_+$, let $B(\la)$ be the crystal associated to an irreducible highest weight $U_q(\g)$-module with highest weight $\la$.
For $w\in W$, the Demazure crystal $B_w(\la)$ is the crystal associated to a $U_q(\mf b)$-module generated by an extremal weight vector $v_{w\la}$ of weight $w\la$, where $U_q(\mf b)$ is the subalgebra of $U_q(\g)$ corresponding to the Borel subalgebra ${\mf b}$ of $\g$ \cite{Kas93}.
The Demazure crystals $B_w(\la)$ for $w\in W$ are characterized as follows:
\begin{equation}\label{eq:Demazure-1}
\begin{split}
B_{e}(\la)&=\{\,v_\la\,\},\ \
B_{s_iw}(\la)=\bigcup_{m\geq 0}\tf_i^m B_w(\la),
\quad\text{for $i\in I$ with $\ell(s_iw)>\ell(w)$},
\end{split}
\end{equation}
where $e$ is the identity in $W$, $v_\la$ is the highest weight element in $B(\la)$, and $\ell(w)$ is the length of $w\in W$.
If $w=s_{i_1}\cdots s_{i_\ell}\in W$ is a reduced expression, then we have
\begin{equation}\label{eq:Demazure-2}
B_w(\la)=\left\{\,\tf_{i_1}^{m_1}\cdots\tf_{i_\ell}^{m_\ell}v_\la \,\Big\vert\,m_1,\ldots,m_\ell\in \Z_+\,\right\}.
\end{equation}
We also have $B_w(\la) \subset B_{w'}(\la)$ if and only if $w\leq w'$
for $w,w'\in W$, where $<$ is a Bruhat order on $W$. Put
\begin{equation*}
\widehat{B}_w(\la)=B_w(\la)\, \setminus \bigcup_{\substack{w'< w \\ w'\la\neq w\la}}B_{w'}(\la).
\end{equation*}
Then we have the following decomposition
\begin{equation}\label{eq:Schubert decomp}
B(\la)=\bigsqcup_{w\in W/W_\la}\widehat{B}_w(\la),
\end{equation}
where $W_\la=\{\,w\in W\,|\,w\la=\la\,\}$.

Suppose that $\g$ is of finite type, and hence $B(\la)$ is a finite crystal for $\la\in P_+$. Let $w_0$ be the longest element in $W$. Similarly, for $w\in W$, the opposite Demazure crystal $B^w(\la)$ is the crystal associated to a $U_q(\mf b^-)$-module generated by an extremal weight vector $v_{w\la}$ of weight $w\la$, where ${\mf b}^-$ is the negative Borel subalgebra of $\g$.
As in \eqref{eq:Demazure-1}, it can be characterized by
\begin{equation}\label{eq:Demazure-2}
\begin{split}
B^{w_0}(\la)&=\{\,v_{w_0\la}\,\},\ \
B^{s_iw}(\la)=\bigcup_{m\geq 0}\te_i^m B^w(\la),
\ \text{for $i\in I$ with $\ell(s_iw)<\ell(w)$}.
\end{split}
\end{equation}
For $\la\in \cP_n$ and $w\in W$, let
\begin{equation*}
\widehat{B}^w(\la)=B^w(\la)\, \setminus \bigcup_{\substack{w'> w \\ w'\la\neq w\la}}B^{w'}(\la).
\end{equation*}
Then we also have the following decomposition
\begin{equation}\label{eq:Schubert decomp-2}
B(\la)=\bigsqcup_{w\in W/W_\la}\widehat{B}^w(\la).
\end{equation}

\subsection{Path model}\label{sec:path model}

Let us briefly recall the Littelmann path model \cite{Li95}. Let $\la\in P_+$ be given.
For $\mu,\nu\in W\la$, we define
%$\mu\leq \nu$
$\nu \geq \mu$
%$\nu \geq \mu$
if there exists a sequence of weights $\nu_0=\nu,\ldots,\nu_s=\mu$ and
a sequence of positive real roots $\beta_1,\ldots,\beta_s$ such that
$\langle \nu_{k-1},h_{\beta_k}\rangle<0$ and $\nu_k=s_{\beta_k}(\nu_{k-1})$,
where $h_{\beta_k}$ and $s_{\beta_k}$ are the coroot and reflection corresponding to $\beta_k$
for $1\leq k\leq s$.
We define
%${\rm dist}(\mu,\nu)$
${\rm dist}(\nu,\mu)$
%{\rm dist}(\nu,\mu)
to be the maximal length $s$ of such sequences.
For $a\in \mathbb{Q}$ and $\mu,\nu\in W\la$ with
%$\mu\leq \nu$,
$\nu \geq \mu,$
%$\nu \geq \mu,$
an $a$-chain for
%$(\mu,\nu)$
$(\nu,\mu)$
%$(\nu,\mu)$
is a sequence
%{\color{blue} $\mu=\la_0 < \cdots  < \la_s=\nu$}
$\nu = \lambda_0 > \lambda_1 > \cdots > \lambda_s = \mu$
%$\nu = \lambda_0 > \lambda_1 > \cdots > \lambda_s = \mu$
of weights in $W\la$ such that
${\rm dist}(\la_{k-1},\la_{k})=1$,
$a\langle \la_{k-1},h_{\beta_k}\rangle\in \Z$,
and $\la_{k}=s_{\beta_k}(\la_{k-1})$ for some positive real root $\beta_k$, for each $1\leq k\leq s$.

Then a Lakshmibai-Seshadri (LS) path of class $\la$ is a pair $\pi=(\underline{\nu};\underline{a})$ of sequences $\underline{\nu} : \nu_0>\cdots>\nu_s$ of weights in $W\la$ and $0=a_0<\cdots<a_s=1$ of rational numbers such that there exists an $a_k$-chain for $(\nu_{k-1},\nu_{k})$ for each $1\leq k\leq s$. Let $\h^*_{\mathbb{R}}=\mathbb{R}\otimes_{\Z} P$. We may regard $\pi$ as a piecewise linear function $\pi : [0,1]\longrightarrow \h^\ast_{\mathbb{R}}$ such that
\begin{equation}\label{eq:LS path}
\pi(t) =\sum_{k=1}^{i-1}(a_k-a_{k-1})\nu_k + (t-a_{i-1})\nu_i,
\end{equation}
for $a_{i-1}\leq t\leq a_i$. We denote by $\B(\la)$ the set of all LS paths of class $\la$.
Let $\pi_\la$ be given by $\pi_\la(t)=t\la$ for $t\in [0,1]$. It is clear that $\pi_\la\in \B(\la)$.

Let $\pi\in \B(\la)$ be given. For $i\in I$, let $h={\rm min}\{\,\langle \pi(t),h_i \rangle\,|\,t\in [0,1]\,\}$. We have $h\in -\Z_+$. If $h=0$, we define $\te_i \pi={\bf 0}$, where ${\bf 0}$ is a formal symbol.
If $h\leq -1$, let $t_1$ be the smallest one such that $\langle \pi(t_1),h_i \rangle =h$, and $t_0\in [0,t_1]$ the largest one such that $\langle \pi(t_0),h_i \rangle =h+1$.
Then we define $\te_i \pi : [0,1]\longrightarrow \h^\ast_\mathbb{R}$ such that
\begin{equation*}
(\te_i \pi)(t)=
\begin{cases}
\pi(t), & \text{for $0\leq t\leq t_0$},\\
\pi(t)-\langle \pi(t)-\pi(t_0),h_i\rangle\alpha_i, & \text{for $t_0\leq t\leq t_1$},\\
\pi(t)+\alpha_i, & \text{for $t_1\leq t\leq 1$}.
\end{cases}
\end{equation*}
We define $\tf_i\pi$ in a similar way. If $\pi(1)-h<1$, we define $\tf_i\pi={\bf 0}$. If $\pi(1)-h\geq 1$,
let $t_1$ be the smallest one such that $\langle \pi(t_1),h_i \rangle =h+1$, and $t_0\in [0,t_1]$ the largest one such that $\langle \pi(t_0),h_i \rangle =h$.
Then we define $\tf_i \pi : [0,1]\longrightarrow \h^\ast_\mathbb{R}$ such that
\begin{equation*}
(\tf_i \pi)(t)=
\begin{cases}
\pi(t), & \text{for $0\leq t\leq t_0$},\\
\pi(t)-\langle \pi(t)-\pi(t_0),h_i\rangle\alpha_i, & \text{for $t_0\leq t\leq t_1$},\\
\pi(t)-\alpha_i, & \text{for $t_1\leq t\leq 1$}.
\end{cases}
\end{equation*}
We put ${\rm wt}(\pi)=\pi(1)$,
$\varepsilon_i(\pi)=\max\{\,m \in \Z_+ \,|\,\te_i^m\pi\neq {\bf 0}\,\}$, and
$\varphi_i(\pi)=\max\{\,m\in\Z_+\,|\,\tf_i^m\pi\neq {\bf 0}\,\}$
for $\pi\in \B(\la)$ and $i\in I$.

\begin{thm}[\cite{Li95}]
For $\la\in P_+$, $\B(\la)$ is a crystal with respect to ${\rm wt}$, $\varepsilon_i$, $\varphi_i$, and $\te_i$, $\tf_i$ for $i\in I$ defined above. Moreover,
%$\B(\la)=\{\,\tf_{i_1}\ldots \tf_{i_r}\pi_\la\,|\,r\geq 0,i_1,\ldots,i_r\in I\,\}$,
\begin{displaymath}
\B(\la)=\left\{\,\tf_{i_1}^{m_1}\ldots \tf_{i_k}^{m_k}\pi_\la \,|\, k\geq 0,\, i_1,\ldots,i_k\in I\,,\,
m_1,\ldots,m_k \in \mathbb{Z}_+ \,\right\},
\end{displaymath}
%\begin{displaymath}
%\B(\la)=\left\{\,\tf_{i_1}^{m_1}\ldots \tf_{i_r}^{m_r}\pi_\la\,|\,r\geq 0,i_1,\ldots,i_r\in I\,,\,
%m_1,\ldots,m_r \in \mathbb{Z}_+ \,\right\},
%\end{displaymath}
where $\pi_\la$ is given by $\pi_\la(t)=t\la$ for $t\in [0,1]$.
\end{thm}

\begin{thm}[\cite{J,Kas96}]\label{thm:pathiso}
For $\la\in P_+$, there is a unique isomorphism $\psi_\la$ of crystals from $B(\la)$ to $\B(\la)$ such that $\psi_\la(v_\la)=\pi_\la$.
\end{thm}

For
%$\pi=(\nu_0,\ldots,\nu_s;a_0,\ldots,a_1)\in \B(\la)$,
$\pi=(\nu_0,\ldots,\nu_s;a_0,\ldots,a_s) \in \B(\la)$,
%$\pi=(\nu_0,\ldots,\nu_s;a_0,\ldots,a_s) \in \B(\la)$,
let $\iota(\pi)=\nu_0$ and $\tau(\pi)=\nu_s$.
If we put
\begin{equation}\label{eq:Demazure cell}
%\begin{split}
%\B_w(\la)&=\{\,\pi\,|\,\pi\in \B(\la),\ \iota(\pi)\leq w\la \,\},\ \
%\widehat{\B}_w(\la)=\{\,\pi\,|\,\pi\in \B(\la),\ \iota(\pi)= w\la \,\},
%\end{split}
\begin{split}
\B_w(\la)&=\left\{\,\pi\in \B(\la) \,\vert \, \iota(\pi)\leq w\la \,\right\},\ \
\widehat{\B}_w(\la)=\left\{\,\pi\in \B(\la) \, \vert \,  \iota(\pi)= w\la \,\right\},
\end{split}
\end{equation}
for $w\in W$, then we have by induction on $\ell(w)$ for $w\in W$ that $\psi_\la(B_w(\la))\subset \B_w(\la)$ and hence by comparing the number of elements of each weight that
\begin{equation}\label{eq:Demazure_crystal_path}
\B_w(\la)=\psi_\la(B_w(\la)),\quad \widehat{\B}_w(\la)=\psi_\la(\widehat{B}_w(\la)).
\end{equation}
Similarly, if $\g$ is of finite type and put
\begin{equation}\label{eq:opposite Demazure cell}
%\begin{split}
%\B^w(\la)&=\{\,\pi\,|\,\pi\in \B(\la),\ \tau(\pi)\geq w\la \,\},\ \
%\widehat{\B}^w(\la)=\{\,\pi\,|\,\pi\in \B(\la),\ \tau(\pi)= w\la \,\}.
%\end{split}
\begin{split}
\B^w(\la)&=\{\, \pi\in \B(\la) \, \vert \, \tau(\pi)\geq w\la \,\},\ \
\widehat{\B}^w(\la)=\{\,\pi\in \B(\la) \, \vert \, \tau(\pi)= w\la \,\}.
\end{split}
\end{equation}
for $w\in W$, then we have
\begin{equation}\label{eq:Demazure_crystal_path}
\B^w(\la)=\psi_\la(B^w(\la)),\quad
\widehat{\B}^w(\la)=\psi_\la(\widehat{B}^w(\la)).
\end{equation}

\section{Non-symmetric Cauchy identity}\label{sec:Cauchy}
Fix a positive integer $n\geq 2$.
In this section, we assume that $\g=\gl_n$ which is spanned by the elementary matrices $e_{ij}$ for $1\leq i,j\leq n$. We have $P^\vee=\bigoplus_{i=1}^n\Z e_{ii}$, $P={\rm Hom}_\Z(P^\vee,\Z)=\bigoplus_{i=1}^n\Z\epsilon_i$, $\Pi^\vee=\{\,h_i:=e_{ii}-e_{i+1\,i+1}\,|\,i\in I\,\}$, and $\Pi=\{\,\alpha_i:=\epsilon_i-\epsilon_{i+1}\,|\,i\in I\,\}$, where $\langle \epsilon_i,e_{jj} \rangle =\delta_{ij}$ for $1\leq i,j\leq n$ and $I=\{1,\ldots, n-1\}$.

\subsection{Bicrystals and diagonal actions}%Review on SST, RSK and bicrystals
Let $\cP_n$ be the set of partitions $\la=(\la_1,\ldots,\la_n)$ of length $\leq n$.
We identify $\la=(\la_1,\ldots,\la_n)\in \cP_n$ with $\la_1\epsilon_1+\cdots+\la_n\epsilon_n\in P_+$.
Let $[n]=\{\,1,\cdots,n\,\}$ and let $\W$ be the set of finite words with letters in $[n]$.
We regard $[n]$ as a crystal $B(\epsilon_1)$,
where ${\rm wt}(k)=\epsilon_k$ for $k\in[n]$.
Identifying ${\bf a}=a_1\ldots a_r\in \W$ with $a_1\otimes \cdots \otimes a_r$, we regard $\W$ as a crystal. The connected component of ${\bf a}$ in $\W$ is isomorphic to $B(\la)$ for some $\la\in \cP_n$. So there exists a unique $b\in B(\la)$ such that $b\equiv {\bf a}$, which we denote by $b_{\bf a}$.

Let
$$\M=\left\{\,M=(m_{ij})_{i,j\in [n]}\,\,\vert\,\, m_{ij}\in\Z_+\,\right\},$$ and
let
$\I$ be the set of biwords $(\ba,\bb)\in \W\times\W$ such that (1) $\ba=a_1\cdots a_r$ and
$\bb=b_1\cdots b_r$  for some $r\geq 0$, (2)   $(a_1,b_1)\leq \cdots
\leq (a_r,b_r)$, where
$(a,b)< (c,d)$ if and only if $(b<d)$ or ($b=d$ and $a>c$), for $a,b,c,d\in [n]$. Then we have
a bijection from $\I$ to $\M$, where
$(\ba,\bb)$ is mapped to $M(\ba,\bb)=(m_{ij})$ with
$m_{ij}=\left|\{\,k\,|\,(a_k,b_k)=(i,j) \,\}\right|$.

For $i\in I$ and $M\in\M$ with $M=M({\bf a},{\bf b})$, we define
\begin{equation}\label{eq:bicrystal-1}
\td{x}_iM=
\begin{cases}
M(\td{x}_i{\bf a},{\bf b}),& \text{if $\td{x}_i{\bf a}\neq {\bf 0}$},\\
{\bf 0}, & \text{if $\td{x}_i{\bf a} = {\bf 0}$},
\end{cases}\quad (x=e,f)
\end{equation}
and put ${\rm wt}(M)={\rm wt}({\bf a})$,
$\varepsilon_i(M)=\varepsilon_i({\bf a})$,
$\varphi_i(M)=\varphi_i({\bf a})$.
Then $\M$ is a crystal with respect to
${\rm wt}$, $\varepsilon_i$, $\varphi_i$, $\te_i$, $\tf_i$, for $i\in I$.
On the other hand, if we define
\begin{equation}\label{eq:bicrystal-2}
\td{x}_i^\sharp M = \left(\td{x}_iM^t \right)^t \quad (x=e,f),
\end{equation}
for $i\in I$ and $M\in \M$, and
${\rm wt}^\sharp(M)={\rm wt}(M^t)$,
%$\varepsilon^\sharp_i(M)=\varepsilon_i(M)$,
%$\varphi^\sharp_i(M)=\varphi_i(M)$
$\varepsilon^\sharp_i(M)=\varepsilon_i(M^t)$,
$\varphi^\sharp_i(M)=\varphi_i(M^t)$
%$\varepsilon^\sharp_i(M)=\varepsilon_i(M^t)$,
%$\varphi^\sharp_i(M)=\varphi_i(M^t)$
($M^t$ denotes the transpose of $M$),
then $\M$ is also a crystal with respect to
${\rm wt}^\sharp$, $\varepsilon^\sharp_i$, $\varphi^\sharp_i$, $\te^\sharp_i$, $\tf^\sharp_i$ for $i\in I$. The following is due to \cite{La} (our presentation here is based on \cite{K07}).

\begin{prop}\label{prop:bicrystal}
Under the above hypothesis,
\begin{itemize}
\item[(1)] $\M$ is a $(\gl_n,\gl_n)$-bicrystal with respect to \eqref{eq:bicrystal-1} and \eqref{eq:bicrystal-2}, that is, $\td{x}_i$ and $\td{y}_j^\sharp$ commute with each other on $\M$ for $x, y\in \{e,f\}$ and $i, j\in I$,

\item[(2)] if $C_\la$ is the connected component of $M_\la:={\rm diag}(\la_1,\ldots,\la_n)$ in $\M$ for $\la\in \cP_n$ as a $(\gl_n,\gl_n)$-bicrystal, then
$$\M=\bigsqcup_{\la\in\cP_n}C_\la.$$

\end{itemize}
\end{prop}

Since $C_\la \cong B(\la)\times B(\la)$ where $\td{x}_i$ and $\td{y}_j^\sharp$ act on the first and second component in $B(\la)\times B(\la)$ for $x, y\in \{e,f\}$ and $i,j\in I$, respectively,
we have an isomorphism of $(\gl_n,\gl_n)$-bicrystals
\begin{equation}\label{eq:RSK}
\begin{split}
\xymatrixcolsep{3pc}\xymatrixrowsep{0.5pc}\xymatrix{
\kappa : \M \ \ar@{->}[r]  & \ \
\displaystyle\bigsqcup_{\la\in\cP_n} B(\la)\times B(\la)  \\
 \quad M \ \ar@{|->}[r]  & \ \ \ \ (b_{\bf a}, b_{\bf c})}
\end{split}
\end{equation}
for $M\in \M$ with $M=M({\bf a},{\bf b})$ and $M^t=M({\bf c}, {\bf d})$. Recall that if we identify $B(\la)$ with the set of semistandard Young tableaux of shape $\la$ with entries in $[n]$, then $\kappa$ is the usual RSK correspondence (cf. \cite{Ful}).

%\subsection{Diagonal action}%State and prove the non-symmetric Cauchy identity
Let ${\mc N}=\bigsqcup_{\la\in\cP_n}B(\la)\times B(\la)$. We define another crystal structure on ${\mc N}$ by regarding $B(\la)\times B(\la)$ as $B(\la)\otimes B(\la)$ as a $\gl_n$-crystal, and denote the associated Kashiwara operators on $\mc N$ by $\td{\tt e}_i$ and $\td{\tt f}_i$ for $i\in I$.
We also define the operators $\td{\tt e}_i$ and $\td{\tt f}_i$ on $\M$ for $i\in I$, which are induced from $\kappa^{-1}$, and hence a crystal structure on $\M$.
By Proposition \ref{prop:bicrystal}, $\kappa$ is an isomorphism of crystals with respect to $\td{\tt e}_i$ and $\td{\tt f}_i$ for $i\in I$.
In fact, if $M\in \M$ is given with $M=M({\bf a},{\bf b})$ and $M^t=M({\bf c}, {\bf d})$,
then
\begin{equation}\label{eq:new quiv}
M\equiv {\bf a}\otimes {\bf c}
\end{equation}
with respect to this crystal structure.
To avoid confusion with the bicrystal structures on $\M$ and $\mc N$, let us denote by $\tt{M}$ and $\tt{N}$ the sets $\M$ and $\mc N$, which are crystals with respect to $\td{\tt e}_i$ and $\td{\tt f}_i$ for $i\in I$.

\subsection{Non-symmetric Cauchy identity}
Let
\begin{equation*}
%\M^{low}=\{\,M\,|\,M=(m_{ij})\in \M,\ m_{ij}=0\ (i<j)\,\}
\M^{low}=\left\{\, M=(m_{ij})\in \M \, \vert \, m_{ij}=0\ (i<j)\,\right\}
%\M^{low}=\left\{\, M=(m_{ij})\in \M \, \vert \, m_{ij}=0\ (i<j)\,\right\}
\end{equation*}
be the set of lower triangular matrices in $\M$.
For $\la\in\cP_n$, let $\tt{C}_\la$ be the connected component of $M_\la$ in $\tt{M}$. We have $M_\la\equiv v_\la\otimes v_\la$ by \eqref{eq:new quiv}, and hence ${\tt C}_\la\cong B(2\la)$.

\begin{lem}\label{lem:decomp of M_low}
Under the above hypothesis,
\begin{itemize}
\item[(1)] $\M^{low}\cup\{{\bf 0}\}$ is invariant under $\td{\tt e}_i$ and $\td{\tt f}_i$ for $i \in I$, where we assume that $\td{\tt e}_i{\bf 0}=\td{\tt f}_i{\bf 0}={\bf 0}$,

\item[(2)] $\M^{low}$ is a subcrystal of $\tt M$, which decomposes as follows:
\begin{equation*}
\M^{low}=\bigsqcup_{\la\in \cP_n}{\tt C}_\la.
\end{equation*}
\end{itemize}
\end{lem}
\pf
(1)
Suppose that there exists $M=(m_{ij})\in \M^{low}$ such that
$\td{\tt e}_iM\neq {\bf 0}$ is not lower triangular for some $i\in I$.
So, if $m_{ii}=x$, $m_{i+1 i}=y$, and $m_{i+1\,i+1}=z$, then we have
\begin{equation}\label{eq:aux-1}
\begin{tikzpicture}[scale=1,baseline=0mm]
\def \adt{160pt}
%%%%%%%%%%%%%%%%%%%%%%%%%%5
\coordinate (A1) at (0,0);
\coordinate (B1) at (\adt,0);

\matrix (A) [matrix of math nodes,
                   column sep=0mm,
                   row sep=-1mm,
                   left delimiter =(,
                   right delimiter =),
                   nodes={
                           line width=1pt,
                          anchor=south,
                          minimum height=2mm}
                  ] at (A1)
{
  {\small \ddots} & \empty & \empty &  \\
   & x & 0 &  \\
  & y & z &  \\
  & \empty & \empty & {\small \ddots} \\
};

\matrix (B) [matrix of math nodes,
                   column sep=-1mm,
                   row sep=-1mm,
                   left delimiter =(,
                   right delimiter =),
                   nodes={
                           line width=1pt,
                          anchor=south,
                          minimum height=2mm}
                  ] at (B1)
{
  {\small \ddots} & \empty & \empty & \empty \\
  \empty & x & 1 & \empty \\
 \empty & y & {\small z-1} & \empty \\
   \empty & \empty & \empty & {\small \ddots} \\
};
\draw[|->] ($(A.east)+(\adt*0.12,0)$) -- ($(B.west)+(-\adt*0.12,0)$) node[above,midway] {$\td{\tt e}_i$};
\end{tikzpicture}
\end{equation}
Equivalently, if $M=M({\bf a},{\bf b})$ and $M^t=M({\bf c}, {\bf d})$ with
\begin{equation*}
\begin{split}
&{\bf a} =  \cdots \underbrace{i+1 \cdots i+1}_{y}  \underbrace{i \cdots i}_{x} \underbrace{i+1 \cdots i+1}_{z}, \\
&{\bf c} =  \underbrace{i \cdots i}_{x} \underbrace{i+1 \cdots i+1}_{z} \underbrace{i \cdots i}_{y}\cdots ,
\end{split}
\end{equation*}
(here we ignore the letters other than $i$ and $i+1$ in $\ba$ and $\bf c$) then $M\equiv {\bf a}\otimes {\bf c}$ and $\widetilde{\tt e}_iM\equiv (\widetilde{\tt e}_i{\bf a})\otimes {\bf c}$, which implies that $x \geq z$ by tensor product rule of crystals (cf.~\cite[Section 2]{KashNaka}).
On the other hand, $\widetilde{\tt e}_i{\bf a}\neq {\bf 0}$ and \eqref{eq:aux-1} implies that $x < z$, which is a contradiction. Hence $\M^{low}\cup\{{\bf 0}\}$ is invariant under $\td{\tt e}_i$ for $i\in I$. The proof for $\td{\tt f}_i$ is similar.

(2) For $\la\in \cP_n$, $M_\la\in \M^{low}$, and hence ${\tt C}_\la\subset \M^{low}$ by (1).  Conversely, let $M=(m_{ij})\in \M^{low}$ be given.
Since $\M^{low}\cup\{{\bf 0}\}$ is invariant under $\td{\tt e}_i$ for $i\in I$, we may assume that $\widetilde{\tt e}_iM={\bf 0}$ for all $i\in I$.
Suppose that $M$ is not diagonal. Then there exists a pair $(i_0,j_0)$ such that $i_0>j_0$ and $m_{i_0j_0}>0$. We may assume that $i_0-j_0=k>0$ is maximal and furthermore $i_0>1$ is the smallest one such $i_0-j_0=k$. By the choice of $(i_0,j_0)$, we have either $j_0=1$ or
\begin{equation*}
m_{i_0 j}=m_{i_0-1\,j}=0\quad (1\leq j\leq j_0-1).
\end{equation*}
\begin{displaymath}
\begin{tikzpicture}[dot/.style={inner sep=0pt,
                                                  minimum size=3pt,
                                                  fill=black,
                                                  circle
                                                  }]
\def \adt{2mm}
%%%%%%%%%%%%%%%%%%%%%%%%%%5
\coordinate (A1) at (0,0);

\matrix (A) [matrix of math nodes,
                   column sep=3mm,
                   row sep=3mm,
                   left delimiter =(,
                   right delimiter =),
                   nodes={
                          minimum width=6mm,
                          anchor=center,
                          minimum height=6mm}
                  ] at (A1)
{
  |[dot]| &\empty  & \empty & \empty & \empty & \empty \\
  \empty & |[dot]| & \empty & \empty & \empty & \empty \\
  \empty & \empty & |[dot]| & \empty & \empty & \empty\\
  \empty & |[dot]| & |[dot]| & |[dot]| & \empty & \empty \\
  \empty & \empty & \empty & |[dot]| & |[dot]|  & \empty \\
  \empty & \empty  & \empty & |[dot]|  & |[dot]|  &  |[dot]|\\
};
\draw[<-] ($(A-4-2) +(-\adt*1,-\adt*0)$) to [out=180,in=45] ($(A-4-2) +(-\adt*8,-\adt*3)$) node[below] {\footnotesize $m_{i_0 j_0}$};
%\draw[-,color=red] ($(A-4-2) +(-\adt*3,-\adt)$) |-  ($(A-3-3) +(\adt*1,\adt*1)$)
%|- ($(A-4-4) +(\adt*1,\adt*1)$) |- ($(A-6-4) +(-\adt*3,-\adt*1)$) |- ($(A-4-2) +(-\adt*3,-\adt)$);
\end{tikzpicture}
\end{displaymath}
This implies that
%$\widetilde{\tt e}_{i_0}M\neq {\bf 0}$,
$\widetilde{\tt e}_{i_0-1}M\neq {\bf 0}$,
%$\widetilde{\tt e}_{i_0-1}M\neq {\bf 0}$,
 which is a contradiction. Hence $M={\rm diag}(\la_1,\ldots,\la_n)$ for some $\la_i\in \Z_+$.
Moreover if $\la_i<\la_{i+1}$ for some $i\in I$, then
%$\widetilde{\tt e}_{i_0}M\neq {\bf 0}$,
$\widetilde{\tt e}_{i}M\neq {\bf 0}$,
%$\widetilde{\tt e}_{i}M\neq {\bf 0}$,
which is also a contradiction. Hence $M=M_\la\in {\tt C}_\la$ with $\la=(\la_1,\ldots,\la_n)\in\cP_n$. This proves (2).
\qed\vskip 2mm

The following is the main theorem in this paper.

\begin{thm}\label{thm:main}
The map $\kappa$ in \eqref{eq:RSK} when restricted to $\M^{low}$ gives a bijection
\begin{equation*}
\xymatrixcolsep{3pc}\xymatrixrowsep{0pc}\xymatrix{
\kappa : \M^{low} \ \ar@{->}[r]  & \ \
\displaystyle \bigsqcup_{\la\in\cP_n}\bigsqcup_{w\in W/W_\la}
B^w(\la)\times \widehat{B}_w(\la)}.
\end{equation*}
\end{thm}
\pf
For $\la\in\cP_n$,
let $N_\la=(v_\la,v_\la)\in \mc N$ and let ${\tt D}_\la$ be the connected component of $N_\la$ in $\texttt N$, which is isomorphic to $B(2\la)$.
By~\eqref{eq:RSK},
%Proposition~\ref{prop:bicrystal},
$\kappa(M_\la)=N_\la$, and hence $\kappa|_{{\tt C}_\la} : {\tt C}_\la \longrightarrow {\tt D}_\la$ is an isomorphism.
So
by Lemma \ref{lem:decomp of M_low},
it suffices to show that
\begin{equation}\label{eq:decomp of D_la}
{\tt D}_\la=\bigsqcup_{w\in W/W_\la}B^w(\la) \times\widehat{B}_w(\la).
\end{equation}
Let $\psi_\la^{-1} : \B(\la) \longrightarrow B(\la)$ be the isomorphism sending $\pi_{\la}$ to $v_\la$.
By \cite[Lemma 2.7]{Li95}, the map sending $\pi\ast \pi'$ to $\pi\otimes\pi'$ is an isomorphism from $\B(\la)\ast\B(\la)$ to $\B(\la)\otimes\B(\la)$, where $\ast$ denotes the concatenation of paths.
So we have a well-defined embedding of crystal
\begin{equation}
\xymatrixcolsep{3pc}\xymatrixrowsep{0pc}\xymatrix{
\theta : \B(2\la) \ \ar@{->}[r]  & \ \ \B(\la)\ast\B(\la)}.
\end{equation}
given by $\theta(\pi)=\pi_1\ast\pi_2$, where
$\pi_1(t)=\pi(t/2)$ and $\pi_2(t)=\pi((t+1)/2)-\pi(1/2)$.
%Note that $\pi_{2\la}=\pi_\la\ast \pi_\la$.
By definition of
the LS path,
we have $\tau(\pi_1)\geq \iota(\pi_2)$, and hence by \eqref{eq:Schubert decomp}
\begin{equation}\label{eq:decomp of B_la}
\theta(\B(2\la))=\bigsqcup_{w\in W/W_\la}\B^w(\la)\ast \widehat{\B}_w(\la).
\end{equation}
Since the map sending $\pi$ to $(\psi_\la^{-1}(\pi_1),\psi_\la^{-1}(\pi_2))$ is an isomorphism
from $\B(2\la)$ to ${\tt D}_\la$, we obtain \eqref{eq:decomp of D_la} from \eqref{eq:decomp of B_la}.
\qed

\begin{cor}
Let
$\M^{up}$
be the set of upper triangular matrices in $\M$.
The map $\kappa$ in \eqref{eq:RSK} when restricted to $\M^{up}$ gives a bijection
\begin{equation*}
\xymatrixcolsep{3pc}\xymatrixrowsep{0pc}\xymatrix{
\kappa : \M^{up} \ \ar@{->}[r]  & \ \
\displaystyle\bigsqcup_{\la\in\cP_n}\bigsqcup_{w\in W/W_\la}
\widehat{B}_w(\la)\times B^w(\la)}.
\end{equation*}
\end{cor}
\pf It follows immediately from the fact that if $\kappa(M)=(b_1,b_2)$ for $M\in \M$, then $\kappa(M^t)=(b_2,b_1)$.
\qed

\begin{cor}\label{cor:main-2}
The map $\kappa$ in \eqref{eq:RSK} when restricted to $\M^{low}$ also gives a bijection
\begin{equation*}
\begin{split}
\xymatrixcolsep{3pc}\xymatrixrowsep{0pc}\xymatrix{
\kappa : \M^{low} \ \ar@{->}[r]  & \ \
\displaystyle\bigsqcup_{\la\in\cP_n}\bigsqcup_{w\in W/W_\la}
\widehat{B}^w(\la)  \times {B}_w(\la)}.
\end{split}
\end{equation*}
\end{cor}
\pf It follows from \eqref{eq:Schubert decomp-2} and the same argument as in Theorem \ref{thm:main}.
\qed\vskip 2mm

Put $z_i=e^{\epsilon_i}\in \mathbb{Q}[P]$ for $1\leq i\leq n$.
For $\la\in \cP_n$ and $w\in W$, let
\begin{equation*}
\begin{split}
&K_{w\la}(z)={\rm ch}B_w(\la),\quad \widehat{K}_{w\la}(z)={\rm ch}\widehat{B}_w(\la),\\
&K^{w\la}(z)={\rm ch}B^w(\la),\quad \widehat{K}^{w\la}(z)={\rm ch}\widehat{B}^w(\la).
\end{split}
\end{equation*}
Recall that $K_{w\la}(z)$ (resp. $K^{w\la}(z)$) is called a Demazure character (resp. opposite Demazure character), while $\widehat{K}_{w\la}(z)$ (resp. $\widehat{K}^{w\la}(z)$) is called a Demazure atom (resp. opposite Demazure atom)
~\cite{Mas09}.

Assume that $x_i, y_j$ are formal commuting variables for $1\leq i,j\leq n$.
From Theorem \ref{thm:main} and Corollary \ref{cor:main-2}, we have the following identities.

\begin{cor}\label{eq:non sym Cauchy id}
We have
\begin{equation*}
\begin{split}
\frac{1}{\prod_{1\leq j\leq i\leq n}(1-x_iy_j)}
&=\sum_{\la\in\cP_n}\sum_{w\in W/W_\la}K^{w\la}(x)\widehat{K}_{w\la}(y)\\
&=\sum_{\la\in\cP_n}\sum_{w\in W/W_\la}\widehat{K}^{w\la}(x){K}_{w\la}(y).
\end{split}
\end{equation*}
\end{cor}

\begin{rem}{\rm
By specializing $y_j=x_j$ for $1\leq j\leq n$, Corollary \ref{eq:non sym Cauchy id} also  recovers the Littlewood identity
\begin{equation}\label{eq:Littlewood identity}
\prod_{1\leq i\leq j\leq n}\frac{1}{(1-x_ix_j)}=\sum_{\la\in\cP_n}s_{2\la}(x_1,\ldots,x_n),
\end{equation}
where $s_\mu(x_1,\ldots,x_n)$ is the Schur polynomial corresponding to $\mu\in \cP_n$.
}
\end{rem}

Now, we recover the expansion of the non-symmetric Cauchy kernel given by Lascoux in \cite[Theorem 6]{La}.
\begin{cor}\label{eq:non sym Cauchy id-2}
We have
\begin{equation*}
\begin{split}
\frac{1}{\prod_{i+ j\leq n+1}(1-x_iy_j)}
&=\sum_{\la\in\cP_n}\sum_{w\in W/W_\la}
{K}_{w\la}(x)\widehat{K}_{ww_0\la}(y).
\end{split}
\end{equation*}
\end{cor}
\pf By replacing $x_i$ with $x_{n-i+1}$ for $1\leq i\leq n$ in Corollary \ref{eq:non sym Cauchy id} and using the fact that $K^{w\la}(x_n,\ldots,x_1)=K_{ww_0\la}(x_1,\ldots,x_n)$, we have
\begin{equation*}
\begin{split}
\frac{1}{\prod_{i+ j\leq n+1}(1-x_iy_j)}
%&=\sum_{\la\in\cP_n}\sum_{w\in W/W_\la}K^{w\la}(x_n,\ldots,x_1)\widehat{K}_{w\la}(y_1,\ldots,y_n)\\
&=\sum_{\la\in\cP_n}\sum_{w\in W/W_\la}
{K}_{ww_0\la}(x)\widehat{K}_{w\la}(y).
\end{split}
\end{equation*}
Finally replacing $\la$ with $w_0\la$ yields the identity.
\qed

\begin{rem}\label{rem:semisimple case}
{\rm
Let $\g$ be a semisimple Lie algebra, and $A_q(\g)$ the quantum coordinate ring associated to $\g$ \cite{Kas93-1}. Then $A_q(\g)$ is a $U_q(\g)$-bimodule and
decomposes as follows:
\begin{equation}\label{eq:Peter-Weyl}
A_q(\g) \cong \bigoplus_{\la\in P_+}V_q(\la)\otimes V_q(\la),
\end{equation}
where $V_q(\la)$ is the irreducible highest weight $U_q(\g)$-module with highest weight $\la$. It also has a crystal base, and its crystal is isomorphic to
$\bigsqcup_{\la\in P_+}B(\la)\times B(\la)$ by \eqref{eq:Peter-Weyl}.

Let $H_q$ be the space of highest weight vectors in $A_q(\g)$ as a $U_q(\g)$-bimodule.
Now, we regard $A_q(\g)$ as a $U_q(\g)$-module where the action is given by
$U_q(\g) \stackrel{\Delta}{\longrightarrow} U_q(\g)\otimes U_q(\g)\longrightarrow {\rm End}(A_q(\g))$ with $\Delta$ a comultiplication. Recall that this corresponds to the diagonal action of $G$ on its coordinate ring where $G$ is the reductive algebraic group associated to $\g$.
 
Let $N_q$ be the $U_q(\g)$-submodule of $A_q(\g)$ generated by $H_q$ and let ${\mc N}_q$ be its crystal.
Since ${\mc N}_q$ is isomorphic to $\bigsqcup_{\la\in P_+}B(2\la)$, the same argument in the proof of \eqref{eq:decomp of D_la} shows that there exists a bijection
\begin{equation*}
\xymatrixcolsep{3pc}\xymatrixrowsep{0pc}\xymatrix{
{\mc N}_q \ \ar@{->}[r]  & \ \
\displaystyle \bigsqcup_{\la\in P_+}\bigsqcup_{w\in W/W_\la}
B^w(\la)\times \widehat{B}_w(\la)}.
\end{equation*}
}
\end{rem}

\section{Continuous crystals and non-symmetric Cauchy kernel}
In this section, we extend the result in the previous section to the case of continuous crystals.

\subsection{Continuous crystals}
Let us review the notion of continuous crystals \cite{BBO,DK,Kas02}.
Let $(A,P^\vee,P,\Pi^\vee,\Pi)$ be the generalized Cartan datum associated to a symmetrizable Kac-Moody algebra as in Section \ref{subsec:crystal}.

A continuous crystal is a set
$B$ together with the maps ${\rm wt} : B \rightarrow P$,
$\varepsilon_i, \varphi_i: B \rightarrow \mathbb{R}\cup\{-\infty\}$ and
$\te^r_i : B \rightarrow B\cup\{{\bf 0}\}$ for $i\in I$ and $r\in\mathbb{R}$ satisfying the following:  for $b\in B$, $i\in I$, and $r\in\mathbb{R}$,
\begin{itemize}
\item[(1)]
$\varphi_i(b) =\langle {\rm wt}(b),h_i \rangle +
\varepsilon_i(b)$,

\item[(2)] $\varepsilon_i(\te^r_i b) = \varepsilon_i(b) - r,\ \varphi_i(\te^r_i b) =
\varphi_i(b) + r,\ {\rm wt}(\te^r_ib)={\rm wt}(b)+r\alpha_i$ if $\te^r_i b \in B$,

%\item[(3)] $\varepsilon_i(\tf_i b) = \varepsilon_i(b) + 1,\ \varphi_i(\tf_i b) = \varphi_i(b) - 1,\ {\rm wt}({\tf_i}b)={\rm wt}(b)-\alpha_i$ if $\tf_i b \in B$,

\item[(3)] $\te^r_i {\bf 0} = {\bf 0}$, $\te^0_i b =b$,

\item[(4)] $\te^{r+s} b =\te^s_i\te^r_i b$ for $s\in \mathbb{R}$, if $\te^r_i b \in B$,

\item[(5)] $\te^r_ib= {\bf 0}$ when $r\neq 0$ and $\varphi_i(b)=-\infty$.
\end{itemize}
For continuous crystals $B_1$ and $B_2$,
a tensor product $B_1\otimes B_2$
is a continuous crystal, which is defined to be $B_1\times B_2$  as a set with elements  denoted by
$b_1\otimes b_2$, where
{\allowdisplaybreaks
\begin{equation*}%\label{eq:tensor product of crystals}
\begin{split}
{\rm wt}(b_1\otimes b_2)&={\rm wt}(b_1)+{\rm wt}(b_2), \\
\varepsilon_i(b_1\otimes b_2)&=
{\rm max}\{\varepsilon_i(b_1),\varepsilon_i(b_2)-\langle {\rm
wt}(b_1),h_i\rangle\}, \\
\varphi_i(b_1\otimes b_2)&= {\rm max}\{\varphi_i(b_1)+\langle {\rm
wt}(b_2),h_i\rangle,\varphi_i(b_2)\},\\
{\te}_i^r(b_1\otimes b_2)&=
\te_i^{r_1}b_1 \otimes \te_i^{r_2} b_2,
\end{split}
\end{equation*}}
for $i\in I$ with
\begin{equation*}
\begin{split}
r_1 &= \max\{r,\varepsilon_i(b_2)-\varphi_i(b_1)\} - \max\{0,\varepsilon_i(b_2)-\varphi_i(b_1)\},\\
r_2 &= \min\{r,\varepsilon_i(b_2)-\varphi_i(b_1)\} + \max\{0,\varphi_i(b_1)-\varepsilon_i(b_2)\}.
\end{split}
\end{equation*}
Here we assume that ${\bf 0}\otimes
b_2=b_1\otimes {\bf 0}={\bf 0}$.

Let ${\mathscr C}$ be the set of continuous maps $\pi : [0,1]\longrightarrow {\mf h}^\ast_{\mathbb R}$. Let $\pi\in \mathscr{C}$ be given. For $i\in I$, put $m_i=\inf_{t\in [0,1]}\{\langle \pi(t), h_i \rangle \}$. For $i\in I$ and $r\in \mathbb{R}$, we define $\te^r_i \pi\in \mathscr{C}$ by
\begin{equation*}
\left(\te^r_i\pi \right)(t) =
\pi (t) -\min\left\{ -r, \inf_{s\in [t,1]}\langle \pi(s), h_i \rangle -m_i \right\}\alpha_i,
\end{equation*}
if $-\langle \pi(1),h_i \rangle +m_i\leq r\leq 0$, and
\begin{equation*}
\left(\te^r_i\pi \right)(t) =
\pi (t) -\min\left\{ 0, -r- m_i+\inf_{s\in [0,t]}\langle \pi(s), h_i \rangle \right\}\alpha_i,
\end{equation*}
if $0\leq r\leq -m_i$, and ${\bf 0}$ otherwise. Put
\begin{equation*}
\begin{split}
&{\rm wt}(\pi)=\pi(1),\quad
\varepsilon_i(\pi)=-m_i,\quad
\varphi_i(\pi)=\langle {\rm wt}(\pi), h_i \rangle + \varepsilon_i(\pi).
\end{split}
\end{equation*}
for $i\in I$. Then $\mathscr{C}$ is a continuous crystal
with respect to ${\rm wt}$, $\varepsilon_i$, $\varphi_i$, $\te^r_i$ for $i\in I$ and $r\in \mathbb{R}$ \cite[Proposition 3.9]{BBO}. In this case, we have
\begin{equation}\label{eq:normality for cont crystals}
\begin{split}
&\varepsilon_i(\pi)=\max\{\,r\in\mathbb{R}_+\,|\,\te_i^r\pi\neq {\bf 0}\,\},\quad
\varphi_i(\pi)=\max\{\,r\in\mathbb{R}_+\,|\,\te_i^{-r}\pi\neq {\bf 0}\,\}.
\end{split}
\end{equation}

Let
%$\ov{C}=\{\,\la\,|\,\la\in\mf{h}^*_{\mathbb{R}},\, \langle \la,h_i \rangle\geq 0\, (i\in I)\,\}$
$\ov{C}=\{\,\la\in\mf{h}^*_{\mathbb{R}}\, \vert \, \langle \la,h_i \rangle\geq 0\  (i\in I)\,\}$
%$\ov{C}=\{\,\la\in\mf{h}^*_{\mathbb{R}}\, \vert \, \langle \la,h_i \rangle\geq 0\, (i\in I)\,\}$
be the closure of the Weyl chamber.
Let $\pi\in \mathscr{C}$ be a continuous path lying in $\ov C$, that is, $\pi(t)\in \ov C$ for all $t\in [0,1]$.
By \eqref{eq:normality for cont crystals}, we have $\te^r_i\pi={\bf 0}$ for $i\in I$ and $r>0$.
In this case, let us put
\begin{equation}\label{eq:cnn comp}
\mathscr{F}(\pi)=
\left\{\,\tf^{r_1}_{i_1}\cdots \tf^{r_k}_{i_k}\pi\,|\,k\geq 0, \,i_1,\ldots i_k\in I,\, r_1,\ldots, r_k\in \mathbb{R}_+\,\right\}.
\end{equation}
Here we denote $\tf^r_i = \te^{-r}_i$ for $i\in I$ and $r\in\mathbb{R}_+$ by convention.

For $\pi, \pi'\in\mathscr{C}$, $\pi\ast \pi'$ denotes the concatenation of $\pi$ and $\pi'$. Then we have the following.
\begin{thm}[Theorem 4.11 in \cite{BBO}]\label{thm:contiuous tensor}
The map
\begin{equation*}
\begin{split}
\xymatrixcolsep{3pc}\xymatrixrowsep{0pc}\xymatrix{
\mathscr{C}\otimes \mathscr{C} \ \ar@{->}[r]  & \ \mathscr{C} \\
\pi\otimes \pi' \ \ar@{|->}[r]  & \ \pi\ast \pi'
}
\end{split}
\end{equation*}
is an isomorphism of continuous crystals.
\end{thm}

Now, let us consider an analogue of Lakshmibai-Seshadri (LS) paths for continuous crystals. We keep the notations in Section \ref{sec:path model}.
Let $\la\in \ov C$ be given.
We define  $\mathscr{B}(\la)$ to be set of the pairs $\pi=(\underline{\nu};\underline{a})$ of sequences $\underline{\nu} : \nu_0>\cdots>\nu_s$ of weights in $W\la$ and $0=a_0<\cdots<a_s=1$ of real numbers.
We regard $\pi$ as a piecewise linear function $\pi : [0,1]\longrightarrow \h^\ast_{\mathbb{R}}$ as in \eqref{eq:LS path}.
We put $\iota(\pi)=\nu_0$ and $\tau(\pi)=\nu_s$. Let $\pi_\la\in \mathscr{B}(\la)$ denote the  path given by $\pi_\la(t)=t\la$ for $t\in [0,1]$.

\begin{prop}\label{prop:continuous LS}
Let $\la\in\ov C$ be given. The set  $\mathscr{B}(\la)\cup\{{\bf 0}\}$ is invariant under $\te^r_i$ for $i\in I$ and $r\in\mathbb{R}$, and hence $\mathscr{B}(\la)$ is a subcrystal of $\mathscr{C}$. Moreover, $\mathscr{B}(\la)=\mathscr{F}(\pi_\la)$.
\end{prop}
\pf Following the arguments in \cite[Section 4]{Li95}, it is straightforward to check that $\mathscr{B}(\la)\cup\{{\bf 0}\}$ is invariant under $\te^r_i$ for $i\in I$ and $r\in\mathbb{R}$. Suppose that $\pi\in \mathscr{B}(\la)$ is given and $\pi\neq \pi_\la$.
Then $\iota(\pi)\neq \la$, and $\langle \iota(\pi),h_i \rangle<0$ for some $i\in I$.
So we have $\varepsilon_i(\pi)>0$, which implies that $\te^{\varepsilon_i(\pi)}_i\pi\in \mathscr{B}(\la)\setminus \{{\bf 0}\}$. By using induction on ${\rm dist}(\iota(\pi),\la)$, we conclude that
there exist
$i_1,\ldots, i_k \in I$
and $r_1,\ldots, r_k\in \mathbb{R}_+$ for some $k\geq 0$ such that
\begin{equation*}
\te^{r_k}_{i_k}\cdots \te^{r_1}_{i_1}\pi=\pi_\la,
\end{equation*}
which is the unique path in $\mathscr{B}(\la)$ lying in $\ov C$.
\qed\vskip 2mm

\begin{prop}\label{prop:conc of cont paths}
For $\la, \mu\in \ov C$, there exists an isomorphism of continuous crystals
\begin{equation*}
\xymatrixcolsep{3pc}\xymatrixrowsep{0pc}\xymatrix{
\mathscr{B}(\la+\mu) \ \ar@{->}[r]  & \ \  \mathscr{F}(\pi_\la\ast\pi_\mu)}
\end{equation*}
sending $\pi_{\la+\mu}$ to $\pi_\la\ast\pi_\mu$.
\end{prop}
\pf Let us give a brief sketch of the proof. Suppose that $\mathscr{C}$ is equipped with a metric $d$ given by $d(\pi,\pi')=\sup_{t\in [0,1]}\{\,||\pi(t)-\pi'(t)||\,\}$,
where $||\cdot ||$ is the usual Euclidean norm on $\mf h^*_\mathbb{R}$.
For $\nu\in \mf h^*_\mathbb{Q}\cap \ov C$, let $\mathscr{B}_\mathbb{Q}(\nu)$ be the subset of $\mathscr{B}(\nu)$ such that $r_1,\ldots,r_k$ in \eqref{eq:cnn comp} are in $\mathbb{Q}$. Then $\mathscr{B}_\mathbb{Q}(\nu)$ is a dense subset of $\mathscr{B}(\nu)$ with respect to $d$. By \cite[Propositions 4.2 and 4.3]{BBO} and \cite[Theorem 5.1]{Li95} one can show that the map
\begin{equation*}
\tf^{r_1}_{i_1}\cdots \tf^{r_k}_{i_k}\pi_{\la+\mu}
\longmapsto
\tf^{r_1}_{i_1}\cdots \tf^{r_k}_{i_k}\pi_\la\ast\pi_\mu
\end{equation*}
for $k\geq 0$, $i_1,\ldots i_k\in I$, and $r_1,\ldots, r_k\in \mathbb{R}_+$,
is a bijection and hence an isomorphism. (Another proof for the case when $\g$ is of finite type can be found in \cite[Theorem 4.11]{BBO}.)
\qed\vskip 2mm

\begin{cor}\label{cor:piecewise linear}
Let $\pi\in \mathscr{C}$ be a piecewise linear path lying in $\ov C$.
Then $\mathscr{F}(\pi)\cong \mathscr{B}(\la)$ as a continuous crystal where $\la=\pi(1)$.
\end{cor}

For $\la\in\ov C$ and $w\in W$, we define $\mathscr{B}_w(\la)$, $\widehat{\mathscr{B}}_w(\la)$,
$\mathscr{B}^w(\la)$, and  $\widehat{\mathscr{B}}^w(\la)$ in the same way as in \eqref{eq:Demazure cell} and \eqref{eq:opposite Demazure cell}.

\subsection{Continuous non-symmetric Cauchy kernel}
From now on, we assume that $\g=\gl_n$ as in Section \ref{sec:Cauchy}.

Let
$$\mathscr{M}=\left\{\,M=(m_{ij})_{i,j\in [n]}\, \vert \, m_{ij}\in\mathbb{R}_+\,\right\}.$$
Let  $M=(m_{ij})_{i,j\in [n]}$ be given. Put
\begin{equation*}
\pi_M = \pi_{M^{(1)}}\ast\cdots \ast \pi_{M^{(n)}} \in \mathscr{C},
\end{equation*}
where $\pi_{M^{(i)}}=(\pi_{m_{n i}\epsilon_n}\ast \pi_{m_{n-1 i}\epsilon_{n-1}}\ast \cdots \ast \pi_{m_{1i}\epsilon_1})$ for $i\in [n]$.
For $i\in I$ and $r\in \mathbb{R}$, we define $\te^r_i M$ to be the unique matrix in $\mathscr{M}$ such that $\te^r_i \pi_M = \pi_{\te^r_i M}$ if $\te^r_i \pi_M \neq {\bf 0}$, and $\te^r_i M={\bf 0}$ if $\te^r_i \pi_M={\bf 0}$.
We put
${\rm wt}(M)={\rm wt}(\pi_M)$,
$\varepsilon_i(M)=\varepsilon_i(\pi_M)$, and
$\varphi_i(M)=\varphi_i(\pi_M)$ for $i\in I$.
Then $\mathscr{M}$ is a continuous crystal with respect to ${\rm wt}$, $\varepsilon_i$, $\varphi_i$, $\te^r_i$ for $i\in I$ and $r\in \mathbb{R}$.

Similarly, given $M\in \mathscr{M}$, we define
$(\te^r_i)^\sharp M = \left(\te^r_iM^t \right)^t$,
${\rm wt}^\sharp(M)={\rm wt}(M^t)$,
%$\varepsilon^\sharp_i(M)=\varepsilon_i(M)$,
%$\varphi^\sharp_i(M)=\varphi_i(M)$ for $i\in I$, $r\in \mathbb{R}$.
$\varepsilon^\sharp_i(M)=\varepsilon_i(M^t)$,
$\varphi^\sharp_i(M)=\varphi_i(M^t)$ for $i\in I$, $r\in \mathbb{R}$.
Then $\mathscr{M}$ is also a continuous crystal with respect to ${\rm wt}^\sharp$,
$\varepsilon^\sharp_i$, $\varphi^\sharp_i$, $(\te^r_i)^\sharp$ for $i\in I$ and $r\in \mathbb{R}$.

As in Proposition \ref{prop:bicrystal}, one can check without difficulty that the operators $\te^r_i$
 and $(\te^s_j)^\sharp$ commute with each other on $\mathscr{M}$ for $i,j\in I$ and $r,s\in \mathbb{R}$, and hence
$\mathscr{M}$ is a continuous $(\gl_n,\gl_n)$-bicrystal.

For $\la=\sum_{i=1}^n\la_i\epsilon_i\in \ov C$, let $M_\la={\rm diag}(\la_1,\ldots,\la_n)$.
By Corollary \ref{cor:piecewise linear}, we have
$\mathscr{F}(M_\la)\cong \mathscr{B}(\la)$ with respect to $\te_i^r$ and
$(\te^r_i)^\sharp$
for $i\in I$ and $r\in \mathbb{R}$ , respectively.

\begin{prop}[cf.\cite{DK}]\label{prop:cont RSK}
There exists an isomorphism of continuous $(\gl_n,\gl_n)$-bicrystals
\begin{equation*}
\xymatrixcolsep{3pc}\xymatrixrowsep{0pc}\xymatrix{
\kappa : \mathscr{M} \ \ar@{->}[r]  &  \ \ \displaystyle\bigsqcup_{\la\in\ov C}\mathscr{B}(\la)\times \mathscr{B}(\la),}
\end{equation*}
such that $\kappa(M_\la)=(\pi_\la,\pi_\la)$ for $\la\in \ov C$,
where $\te^r_i$ and $(\te_i^r)^\sharp$ for $i\in I$ and $r\in \mathbb{R}$ act on the first and second component in $\mathscr{B}(\la)\times \mathscr{B}(\la)$, respectively.
\end{prop}
\pf Let $M\in \mathscr{M}$ be given.
As in Proposition \ref{prop:bicrystal} we can check that there exists $\la\in \ov C$ such that
\begin{equation*}
\te^{r_k}_{i_k}\cdots \te^{r_1}_{i_1}\te^{s_l}_{j_l}\cdots \te^{s_1}_{j_1}M=M_\la
\end{equation*}
for some $i_1,\ldots,i_k$, $j_1,\ldots,j_l\in I$ and $r_1,\ldots,r_k$, $s_1,\ldots,s_l\in \mathbb{R}_+$. This implies that $\mathscr{M}$ is the union of the connected components of $M_\la$ for $\la\in\ov C$ and proves the decomposition of $\mathscr{M}$ as a continuous $(\gl_n,\gl_n)$-bicrystal.
\qed\vskip 2mm

Let $\mathscr{N}=\bigsqcup_{\la\in\ov C}\mathscr{B}(\la)\times \mathscr{B}(\la)$.
We define a continuous crystal structure on $\mathscr{N}$ by identifying $(\pi,\pi')$ with $\pi\otimes \pi'$.
By Proposition \ref{prop:cont RSK}, we can define a continuous crystal structure on $\mathscr{M}$ such that $\kappa : \mathscr{M}\longrightarrow \mathscr{N}$ is an isomorphism of continuous crystals.

Let
\begin{equation*}
\mathscr{M}^{low}=\left\{\, M=(m_{ij})\in \mathscr{M}\, \vert \, m_{ij}=0\ (i<j)\,\right\}.
\end{equation*}
Then we have the following analogue of Theorem \ref{thm:main} for continuous crystals.

\begin{thm}\label{thm:main-2}
The map $\kappa$ when restricted to $\mathscr{M}^{low}$ gives a bijection
\begin{equation*}
\begin{split}
\xymatrixcolsep{3pc}\xymatrixrowsep{0pc}\xymatrix{
\kappa : \mathscr{M}^{low} \ \ar@{->}[r]  & \ \
\displaystyle \bigsqcup_{\la\in\ov C}\bigsqcup_{w\in W/W_\la}
\mathscr{B}^w(\la)\times \widehat{\mathscr{B}}_w(\la)}
\end{split}
\end{equation*}
\end{thm}
\pf The proof is almost identical to that of Theorem \ref{thm:main}. We leave the details to the readers.
\qed

\begin{rem}{\rm
Taking the intersection $\mathscr{M}^{low}\cap \M$ in Theorem \ref{thm:main-2} recovers Theorem \ref{thm:main}.
}
\end{rem}

\end{document}